\documentclass[reqno, oneside, 12pt]{amsart}

\usepackage[letterpaper]{geometry}
\usepackage{array}

\usepackage{enumerate, hyperref,url,amssymb,amsmath,amsthm,amsxtra,mathtools,mathrsfs,calc,nccmath,color}

\usepackage{calc}
\usepackage{graphicx}
\usepackage{caption}
\usepackage{subcaption}
\usepackage{relsize}

\newcommand{\Z}{\mathbb{Z}}
\newcommand{\Q}{\mathbb{Q}}
\newcommand{\R}{\mathbb{R}}

\newcommand{\C}{\mathbb{C}}

\let\temp\phi
\let\phi\varphi
\let\varphi\temp

\renewcommand{\(}{\left(}
\renewcommand{\)}{\right)}

\newcommand{\GL}{\operatorname{GL}}

\renewcommand{\sl}{\big|}
\newcommand{\sk}{\big|_k }

\newcommand{\pMatrix}[4]{\left(\begin{matrix}#1 & #2 \\ #3 & #4\end{matrix}\right)}

\renewcommand{\bar}[1]{\overline{#1}}
\newcommand{\floor}[1]{\left\lfloor #1 \right\rfloor}

\newtheorem{theorem}{Theorem}[section]
\newtheorem{lemma}[theorem]{Lemma}

\theoremstyle{remark}
\newtheorem*{remark}{Remark}

\numberwithin{equation}{section}


\begin{document}


\title{Zeros of period polynomials for symmetric power $L$-functions}

\author{Robert Dicks}
\address{School of Mathematical and Statistical Sciences\\
Clemson University\\
Clemson, SC 29634-0975\\
USA}
\email{rdicks@clemson.edu} 

\author{Hui Xue}
\address{School of Mathematical and Statistical Sciences\\
Clemson University\\
Clemson, SC 29634-0975\\
USA}
\email{huixue@clemson.edu}

 
\begin{abstract}
Suppose that $k$ and $N$ are positive integers.
Let $f$ be a newform on $\Gamma_0(N)$ of weight $k$ with $L$-function $L_f(s)$.
 Previous works have studied the zeros of the period polynomial $r_f(z)$, which is a generating function for the critical values  of $L_f(s)$ and has a functional equation relating $z$ and $-1/Nz$.
 In particular, $r_f(z)$ satisfies a version of the Riemann hypothesis: all of its zeros are on the circle of symmetry
 $\{z \in \C \ : \ |z|=1/\sqrt{N}\}$.

In this paper, for a positive integer $m$, we define a natural analogue of $r_f(z)$ for the $m^{\operatorname{th}}$ symmetric power $L$-function of $f$ when $N$ is squarefree. Our analogue also has a functional equation relating $z$ and $-1/Nz$. We prove the corresponding version of the Riemann hypothesis when $k$ is large enough. Moreover, when $k>2(\operatorname{log}_2(13e^{2\pi}/9)+m)+1$, we prove our result when $N$ is large enough. 
\end{abstract}


\maketitle
\section{Introduction}
For positive integers $k$ 
and $N$,
let $S_k(N)$ denote the complex vector space of cusp forms of weight $k$ on $\Gamma_0(N)$.
 Let 
 $f \in S_k(N)$
 be a newform with Fourier expansion
 $f=\sum_{n \geq 1} a(n)q^{n}$ (see Section $2$ for details).
The $L$-function for $f$, which is given by
\begin{equation}\label{Lfunction}
L_f(s):= \sum^{\infty}_{n=1}\frac{a(n)}{n^s}=\(\prod_{p \mid N}\frac{1}{1-a(p)p^{-s}}\)\prod_{p \nmid N} \frac{1}{1-a(p)p^{-s}+p^{k-1-2s}},
\end{equation}
has an analytic continuation to the complex plane.
Moreover, the completed $L$-function
\[
L^{*}_f(s):= \(\frac{\sqrt{N}}{2 \pi}\)^s\Gamma(s)L_f(s)
\]
satisfies the functional equation
\begin{equation}\label{funceq1}
 L^{*}_f(s)=\epsilon_fL^{*}_f(k-s)
 \end{equation}
for some $\epsilon_f \in \{\pm 1\}$.

The completed $L$-function can be written as 
\begin{equation}\label{perint}
L^{*}_f(s)=N^{s/2}\int^{\infty}_{0}f(iy)y^{s}\frac{dy}{y},
\end{equation}
which means that it arises as a ``period integral" of $f$ (see \cite{Zagier/Kontsevich} for more information on period integrals).
The period polynomial of $f$ is the degree $k-2$ polynomial
\[
r_f(z):= \int^{i\infty}_{0} f(\tau)(\tau-z)^{k-2}d\tau.
\]
Using \eqref{perint} and expanding $(\tau-z)^{k-2}$, we can write $r_f(z)$ as
\begin{equation}\label{perintalt}
r_f(z)=i^{k-1}N^{-\frac{k-1}{2}}\sum^{k-2}_{n=0}{\binom{k-2}{n}}(\sqrt{N}iz)^nL^{*}_f(k-1-n)
\end{equation}
or as
\[
r_f(z)=-\frac{(k-2)!}{(2 \pi i)^{k-1}}\sum^{k-2}_{n=0}\frac{(2 \pi i z)^n}{n!}L_f(k-1-n).
\]
Thus, $r_f(z)$ is a generating function for the critical values $L_f(1)$, $L_f(2)$, $\dots$, $L_f(k-1)$. These values are crucially important in arithmetic geometry and number theory \cite{Ono-Sound1} \cite{Thorner}.
For example, they appear in the context of the Bloch-Kato conjecture for special values of $L$-functions which generalizes the Birch and Swinnerton-Dyer conjecture relating Tate-Shafarevich groups of elliptic curves to analytic properties of modular $L$-functions (see \cite{Bloch-Kato} for more details). 
Moreover, by virtue of the functional equation \eqref{funceq1} and \eqref{perintalt}, we see that
\begin{equation}\label{Perpolyfunc}
r_f(z)=-\epsilon_fi^k(\sqrt{N}z)^{k-2}r_f\(-\frac{1}{Nz}\).\footnote{Here, we have corrected a typo in \cite{Ono-Sound1}; their version of \eqref{Perpolyfunc} had $\frac{k-2}{2}$ instead of $k-2$.}
\end{equation}

 Jin, Ma, Ono, and Soundararajan \cite{Ono-Sound1} showed that $r_f(z)$ satisfies an analogue of the Riemann hypothesis: all of its zeros are located on the circle $\{z \in \C \ : \ |z|=1/\sqrt{N}\}$, the circle of symmetry for the functional equation \eqref{Perpolyfunc}. 
 In this paper, for a positive integer $m$, we define an analogue of $r_f(z)$ for the $m^{\operatorname{th}}$ symmetric power $L$-function and prove the corresponding version of the Riemann hypothesis when $N$ is squarefree and when either $k$ is large enough, or $N$ is large enough assuming that $k>2(\operatorname{log}_2(13e^{2\pi}/9)+m)+1$.

Stating our result requires introducing notation. 
Suppose that $k$ is a positive integer, that $N$ is a squarefree positive integer and that
$f =\sum_{n \geq 1} a(n)q^{n} \in S_k(N)$ is a newform.
 For each prime $p$, we denote by $\alpha_p$ and $\beta_p$ the roots of the polynomial
\[
X^2-a(p)X+p^{k-1}
\]
when $p \nmid N$; otherwise, we define $\alpha_p=a(p)$ and $\beta_p=0$.
Let $m$ be a positive integer.
 The  $m^{\operatorname{th}}$ symmetric power $L$-function of $f$ is defined by
\begin{equation}\label{Euler product}
L_{m,f}(s):= \prod_{p} \prod_{0 \leq i \leq m} (1-\alpha^{m-i}_p \beta^i_p p^{-s})^{-1},
\end{equation}
which converges absolutely for ${\rm Re}(s) > \frac{m(k-1)}{2}+1$ \cite[pg.~$224$]{Xiao}. Note that when $m=1$ we have $L_{1,f}(s)=L_f(s)$.

 Newton and Thorne \cite[Theorem $A$]{Newton-Thorne} showed that $L_{m,f}(s)$ has an analytic continuation to the whole complex plane. 
 Moreover, the completed  $m^{\operatorname{th}}$ symmetric power $L$-function
\[
L_{m,f}^{*}(s):=N^{\frac{ms}{2}}\gamma_m(s)L_{m,f}(s),
\]
satisfies the functional equation
\begin{equation}\label{compfunc}
L_{m,f}^{*}(s)= \epsilon_{m,f}L^{*}_{m,f}(m(k-1)+1-s)
\end{equation}
for some $\epsilon_{m,f} \in \{\pm 1\}$.
Here \cite[pg. $118$]{Zagier}
\begin{equation}\label{Zagier}
\gamma_m(s):= 
\begin{cases} 
 (2 \pi)^{-rs}\displaystyle \prod^{r-1}_{j=0} \Gamma(s-j(k-1))& \text{if } m=2r-1, \\
\pi^{-\frac{s}{2}}\Gamma\(\frac{s}{2}-\floor{\frac{r(k-1)}{2}}\)\gamma_{2r-1}(s) & \text{if } m=2r.
\end{cases}
\end{equation}
Zagier defines $\gamma_m(s)$ only for $N=1$, but it does not depend on $N$ (see e.g. \cite[$(1.6)$]{Xiao}, where the gamma factor differs from ours by a constant multiple). 

The main object of our investigation is 
\begin{equation}\label{symmperpoly}
R_{m,f}(z):= i^{m(k-1)}N^{-\frac{m(k-1)}{2}}\mathlarger{\sum}^{\mathsmaller{m(k-1)-1}}_{n=0} {\binom{m(k-1)-1}{n}} (\sqrt{N}iz)^n L^{*}_{m,f}(m(k-1)-n),
\end{equation}
which is our analogue of $r_f(z)$. In fact, it follows from \eqref{perintalt} that $R_{1,f}(z)=r_f(z)$; moreover
the polynomial $R_{m,f}(z)$ satisfies a functional equation similar to \eqref{Perpolyfunc} (see Lemma~\ref{2.3} below). 
 With this notation, our main result is as follows.
\begin{theorem}\label{main theorem}
Suppose that $k$, $m$ and $N$ are positive integers such that $m>1$ and $N$ is squarefree. Let $f \in S_k(N)$ be a newform. 
Then $R_{m,f}(z)$ has all of its zeros on the circle $\{z \in \C \ : \ |z|=1/\sqrt{N}\}$ if $k$ is large enough. The result also holds for large enough $N$ if $k>2(\operatorname{log}_2(13e^{2\pi}/9)+m)+1$.
\end{theorem}

\begin{remark}
(1) When $m=1$, the result is true for all $k$ and $N$, see \cite[Theorem $1.1$]{Ono-Sound1}.\\
(2) Löbrich, Ma, and Thorner \cite{Thorner} have done related work for special values of motivic $L$-functions. For a pure motive $\mathcal{M}$ over $\Q$ of rank $d \geq 2$, global conductor $N \in \Z^{+}$ and odd weight $w$, they define a polynomial $p_{\mathcal{M}}(z)$ which resembles $r_f(z/i\sqrt{N})$ by using the completed motivic $L$-function $\Lambda(s,\mathcal{M})$ and the Hodge numbers for $\mathcal{M}$. Assuming that the (usual) motivic $L$-function $L(s,\mathcal{M})$ coincides with the $L$-function of an algebraic tempered cuspidal symplectic representation of $\GL_d(\mathbb{A}_{\Q})$, they prove for a large number of cases that the zeros are on the unit circle and are equidistributed. Moreover, their work has applications to odd symmetric power $L$-functions of cuspidal newforms. More precisely, when $f$ is a non-CM newform of squarefree level $N \geq 13$ and weight $k \geq 6$ with integral coefficients whose $L$-function coincides with that of an algebraic tempered cuspidal symplectic representation of $\GL_{n+1}(\mathbb{A}_{\Q})$, then all of the zeros of their polynomial are on the unit circle and they are equidistributed as $n$ or $N$ goes to $\infty$ (they have more complicated results in weights $2$ and $4$).

Our work differs from theirs in that our definition is slightly simpler (the key difference being the lack of an analogue of Hodge numbers in our definition) and we have no restriction on the Fourier coefficients of the forms under our consideration. \\
(3) 
Computations in Sage for $N=1$ show that the zeros of $R_{2,f}(z)$ are very close to the unit circle (if not on it). It is natural to suspect that Theorem~\ref{main theorem} holds in this case and to investigate whether or not it holds for other small values of $k$ and $N$. \\
(4)
The obstruction to proving the result when $N$ is not squarefree is that we do not have an explicit Euler product as in \eqref{Euler product} in this case  (see Lemma~\ref{lastprep} below).
\end{remark}
The main tool that we use in the proof of Theorem~\ref{main theorem} is the following result of Lal\`in and Smyth 
(extending work of \cite{Cohn}) on locating the zeros of ``self-inversive" polynomials (for more background on these polynomials, see \cite[pp. $201-206$]{Marden}).
\begin{theorem}{\cite[Theorem $1$]{LalinSmyth}}\label{ElG,R}
Let $n$ be a positive integer and $h(z) \in \C[z]$ be a  polynomial of degree $n$ with all of its zeros in $\{z \in \C \ : \ |z| \leq 1\}$.
 If $d \geq n$ is a positive integer and $\lambda \in \C$, set
\begin{equation}\label{selfinversive}
P^{\lambda}(z):= z^{d-n}h(z)+\lambda z^{n}\bar{h}(z^{-1}).
\end{equation}
Then, if $|\lambda|=1$ and $P^{\lambda} \neq 0$, the polynomial $P^{\lambda}$ has all of its zeros on the unit circle.
\end{theorem}

The paper is organized as follows.
In Section $2$, we give background on modular forms and the estimates on $L$-functions which we will need. We also prove a functional equation for \eqref{symmperpoly} which we will need in Section $3$. 
In Section $3$, we prove Theorem~\ref{main theorem} for odd symmetric powers (the proofs for the odd and even symmetric powers are slightly different because of the definition of the gamma factor in \eqref{Zagier}). We construct $P_{m,f}(z) \in \R[z]$ such that all of the zeros of $R_{m,f}(z)$ are on the circle $\{z \in \C \ : |z|=1/\sqrt{N}\}$ if and only if $P_{m,f}(z)$ has all of its zeros on the unit circle.
We then construct $Q_{m,f}(z)\in \R[z]$ such that $P_{m,f}(z)$ and $Q_{m,f}(z)$ satisfy the hypotheses of Theorem~\ref{ElG,R}.
That the zeros of $Q_{m,f}(z)$ satisfy $|z| \leq 1$ is a consequence of Rouch\'e's theorem (see Lemma~\ref{mainlemma} below, where we can in fact say that the zeros satisfy $|z|<1$). 
In Section $4$, we sketch the proof of Theorem~\ref{main theorem} for even symmetric powers; the details are similar to those of Section $3$ once we modify the constructions which appear there.

\section{Background}
For background on modular forms, see \cite{Diamond-Shurman}.
Suppose that $k$ and $N$ are positive integers.
  For a function $f(z)$ on the upper half plane  and 
\[
 \gamma =\left(\begin{matrix}a & b \\c & d\end{matrix}\right) \in \GL_{2}^{+}(\Q),
\]
we have the weight $k$ slash operator 
\[
f(z)\sk \gamma := \det(\gamma)^{\frac{k}{2}}(cz+d)^{-k}f\(\frac{az+b}{cz+d}\).
\]
We denote by $S_k(N)$ the space of cusp forms of weight $k$ on $\Gamma_0(N)$. Forms in this space satisfy the transformation law
\[
f \sl_k \gamma=f \ \ \text{ for } \ \ \gamma=\pMatrix{a}{b}{c}{d} \in \Gamma_0(N)
\]
and vanish at the cusps of $\Gamma_0(N)$. 


If $p$ is prime, we denote by $T_p$ the usual Hecke operator on $S_k(N)$ (see \cite[Chapter $5$]{Diamond-Shurman}). As in \cite[Chapter $5$]{Diamond-Shurman}, we use $T_p$ even when $p \mid N$. We say that $f= \sum_{n \geq 1} a(n)q^{n}\in S_k(N)$ is a newform if $a(1)=1$ and it is an eigenform of all the $T_p$ operators.

We require the following estimate for symmetric power $L$-functions (which is essentially a special case of \cite[Lemma $2.2$]{Xiao}).
\begin{lemma}\label{convex}
Suppose that $m$ is a positive integer, that $N$ is a squarefree positive integer and that $k$ is a positive integer. Let $f \in S_{k}(N)$ be a newform. 
We have
\[
L_{m,f}(s) \ll_{m} \operatorname{log}^{m+1}\(k+\Big|s+\frac{m(k-1)}{2}-1\Big|+2\),
\]
uniformly for $s \geq \frac{m(k-1)}{2}$.
\end{lemma}
\begin{remark}
(1) Lemma $2.2$ of \cite{Xiao} has the additional hypothesis that $1 \leq m \leq 4$. However, the remark on page $228$ of \cite{Xiao} mentions that his results apply for the symmetric power $L$-functions of $f$ for which automorphy and cuspidality are known, which is now the case for all $m$ by \cite[Theorem $A$]{Newton-Thorne}.\\
(2) That the implicit constant depends on $m$ alone is missing from the statement of \cite[Lemma $2.2$]{Xiao}, but it follows from its proof.
\end{remark}
We also need a stronger estimate for the proof of Lemma~\ref{mainlemma} below. To describe this, for positive integers $n$ and $w$ with $w \geq 2$, let $\textbf{1}(n)$ denote the constant function whose value is $1$ and let $d_w(n)$ denote the $w^{\operatorname{th}}$ convolution power of $\textbf{1}(n)$. 

\begin{lemma}\label{lastprep}
Suppose that $m$ is a positive integer, that $N$ is a squarefree positive integer and that $k$ is a positive integer with $k \geq 6$. Let $f \in S_k(N)$ be a newform. Then for $s \geq \frac{(m+1){k-1}}{2}$ we have 
\[
|L_{m,f}(s)-1| <\frac{13}{9} \times 2^{m-\frac{k-1}{2}}.
\]
\end{lemma}
\begin{proof}
With minor changes, we follow the proof of the first assertion of \cite[Lemma $2.4$]{Farmer-Conrey-Imamoglu}.
By \cite[$(2.1)$]{Xiao} and a standard convolution identity for the Riemann zeta function, we have
\[
|L_{m,f}(s)-1| \leq \sum^{\infty}_{n=2} \frac{d_{m+1}(n)}{n^{\frac{k-1}{2}}}=\zeta\(\frac{k-1}{2}\)^m-1.
\]
Since $k \geq 6$, we have
$\zeta\(\frac{k-1}{2}\)^m-1< 2^m\(\zeta\(\frac{k-1}{2}\)-1\)$ 
and
\[
\zeta\(\frac{k-1}{2}\)-1=2^{-\frac{k-1}{2}}+\sum^{\infty}_{n=3} n^{-\frac{k-1}{2}} \leq 2^{-\frac{k-1}{2}}+\int^{\infty}_{2}u^{-\frac{k-1}{2}}du \leq \frac{13}{9} \times 2^{-\frac{k-1}{2}}.
\]
The result follows.
\end{proof}
We now describe and prove the functional equation for \eqref{symmperpoly} mentioned in the introduction.
\begin{lemma}\label{2.3}
Suppose that $k$ and $m$ are positive integers, that $N$ is a squarefree positive integer and that
$f \in S_k(N)$ is a newform.
Let $R_{m,f}(z)$ be given by \eqref{symmperpoly}. We have 
\begin{equation}\label{perpolyfunc}
R_{m,f}(-1/Nz)=\epsilon_{m,f}i^{3(m(k-1)-1)}z^{1-m(k-1)}N^{\frac{1-m(k-1)}{2}}R_{m,f}(z).
\end{equation}
\end{lemma}
\begin{proof}
From \eqref{symmperpoly}
\begin{equation*}
R_{m,f}(z)= i^{m(k-1)}N^{-\frac{m(k-1)}{2}}\mathlarger{\sum}^{\mathsmaller{m(k-1)-1}}_{n=0} {\binom{m(k-1)-1}{n}} (\sqrt{N}iz)^n L^{*}_{m,f}(m(k-1)-n),
\end{equation*}
and the functional equation \eqref{compfunc}
\begin{equation*}
L_{m,f}^{*}(s)= \epsilon_{m,f}L^{*}_{m,f}(m(k-1)+1-s),
\end{equation*} we see that 
{\scriptsize{\begin{align*}
 &R_{m,f}(-1/Nz)
\\=&i^{m(k-1)}z^{1-m(k-1)} \sum^{m(k-1)-1}_{n=0}\binom{m(k-1)-1}{n}i^{3n}N^{-\frac{m(k-1)+n}{2}}z^{m(k-1)-1-n}L^{*}_{m,f}(m(k-1)-n)\\
=&i^{m(k-1)}z^{1-m(k-1)}\epsilon_{m,f} \sum^{m(k-1)-1}_{n=0}\binom{m(k-1)-1}{n}i^{3n}N^{-\frac{m(k-1)+n}{2}}z^{m(k-1)-1-n}L^{*}_{m,f}(n+1)
\\
=&i^{m(k-1)}z^{1-m(k-1)}\epsilon_{m,f} 
\sum^{m(k-1)-1}_{n=0}\binom{m(k-1)-1}{n}i^{3(m(k-1)-1-n)}N^{-m(k-1)+\frac{n+1}{2}}z^{n}L^{*}_{m,f}(m(k-1)-n).\\
=&\epsilon_{m,f}i^{3(m(k-1)-1)}z^{1-m(k-1)}N^{\frac{1-m(k-1)}{2}}i^{m(k-1)} 
N^{-\frac{m(k-1)}{2}}\sum^{m(k-1)-1}_{n=0}\binom{m(k-1)-1}{n}(\sqrt{N}iz)^{n}L^{*}_{m,f}(m(k-1)-n),
\end{align*}}}
which concludes the proof.
\end{proof}
\section{Proof of Theorem~\ref{main theorem} for odd symmetric powers}
Unless otherwise stated, we assume that $m > 1$ is an odd integer with $m=2r-1$, that $k$ and $N$ are positive integers with $N$ squarefree and that $f \in S_k(N)$ is a newform.
For each of the summands of $R_{m,f}(z)$ in \eqref{symmperpoly}, we compute that
\begin{multline}
\binom{m(k-1)-1}{n}L_{m,f}^{*}(m(k-1)-n)=\\
\mfrac{(m(k-1)-1)!}{n!}(2 \pi)^{-r(m(k-1)-n)}\(\prod^{r-1}_{j=1}\Gamma((m-j)(k-1)-n)\)N^{\frac{m^2(k-1)-mn}{2}} L_{m,f}(m(k-1)-n).
\end{multline}
With a view towards using Theorem~\ref{ElG,R}, we define
\[
C_{m,f}:=(2 \pi)^{rm(k-1)}(m(k-1)-1)!^{-1}\displaystyle \(\prod^{r-1}_{j=1}\Gamma((m-j)(k-1))\)^{-1}N^{-\frac{m^2(k-1)}{2}-1}
\]
and
\begin{equation}\label{realcoeffpoly}
P_{m,f}(z):=C_{m,f}R_{m,f}(z/i\sqrt{N}).
\end{equation}
The zeros of $R_{m,f}(z)$ are all on the circle $\{z \in \C \ : \ |z|=1/\sqrt{N}\}$ if and only if the zeros of $P_{m,f}(z)$ are on the unit circle.
Moreover, 
by \eqref{perpolyfunc}, we have the functional equation
\begin{equation}\label{realcoefffunc}
P_{m,f}(1/z)=\epsilon_{m,f}z^{1-m(k-1)}P_{m,f}(z).
\end{equation}

If we set
\begin{multline}\label{Qfm}
Q_{m,f}(z):=
\frac{1}{2}C_{m,f}\binom{\mathsmaller{m(k-1)-1}}{\frac{m(k-1)-1}{2}}L_{m,f}^{*}\(\frac{m(k-1)+1}{2}\)+\\
C_{m,f} \mathlarger{\sum}^{\mathsmaller{\frac{m(k-1)-1}{2}-1}}_{n=0} \binom{m(k-1)-1}{n} L_{m,f}^{*}\(m(k-1)-n\)z^{\frac{m(k-1)-1}{2}-n},
\end{multline}
then the functional equations \eqref{compfunc} and \eqref{realcoefffunc} imply that
\begin{equation}\label{useELG,R}
P_{m,f}(z)=z^{\frac{m(k-1)-1}{2}}Q_{m,f}(z)+\epsilon_{m,f}z^{\frac{m(k-1)-1}{2}}Q_{m,f}(1/z).
\end{equation}
Thus, the core of our proof is showing that the zeros of $Q_{m,f}(z)$ satisfy $|z|<1$.
 To this end, we begin with the following lemma. For any positive integer $N$, define 
\begin{multline}
H_{m,f}(z):=
\frac{(2\pi)^{\frac{m(k-1)-1}{2}}}{2(\frac{m(k-1)-1}{2})!}\(\prod^{r-1}_{j=1}\frac{\Gamma(\frac{m(k-1)-1}{2}+1-j(k-1))}{\Gamma((m-j)(k-1))}\)N^{-\frac{m^2(k-1)}{2}-1}+\\
\mathlarger{\sum}^{\frac{m(k-1)-1}{2}-1}_{n=0}\frac{(2\pi)^n}{n!}\(\prod^{r-1}_{j=1}\frac{\Gamma((m-j)(k-1)-n)}{\Gamma((m-j)(k-1))}\)N^{-\frac{n}{2}-1}z^{\frac{m(k-1)-1}{2}-n}.\label{eq:Hmf}
\end{multline}
With this notation, we have the following result.
\begin{lemma}\label{preplemma}
Suppose that $m>1$ is an odd integer with $m=2r-1$ that $k$ and $N$ are positive integers, and that $f \in S_k(N)$ is a newform.
If $k+N$ is sufficiently large (independent of $m$) then all of the zeros of $H_{m,f}(z)$ satisfy $|z|<1$.
\end{lemma}
\begin{proof}
We pick the highest two degree terms of \eqref{eq:Hmf} and define
\begin{equation}\label{truncate}
M_{m,f}(z):=N^{-1}z^{\frac{m(k-1)-1}{2}}+2 \pi\(\prod^{r-1}_{j=1}\frac{\Gamma((m-j)(k-1)-1)}{\Gamma((m-j)(k-1))}\)N^{-\frac{3}{2}}z^{\frac{m(k-1)-1}{2}-1}.
\end{equation}
We will compare the values of $M_{m,f}(z)$ and $H_{m,f}(z)-M_{m,f}(z)$ on the unit circle.
Note that $\Gamma(h)=(h-1)!$ for $h \in \Z^{+}$.
We have
\begin{multline}
|H_{m,f}(z)-M_{m,f}(z)|=
\mathlarger{\mathlarger{\mathlarger{\sl}}} \frac{(2\pi)^{\frac{m(k-1)-1}{2}}}{2(\frac{m(k-1)-1}{2})!}\(\prod^{r-1}_{j=1}\frac{\Gamma(\frac{m(k-1)-1}{2}+1-j(k-1))}{\Gamma((m-j)(k-1))}\)N^{-\frac{m^2(k-1)}{2}-1}+\\
\mathlarger{\sum}^{\frac{m(k-1)-1}{2}-1}_{n=2}\frac{(2\pi)^n}{n!}\(\prod^{r-1}_{j=1}\frac{\Gamma((m-j)(k-1)-n)}{\Gamma((m-j)(k-1))}\)N^{-\frac{n}{2}-1}z^{\frac{m(k-1)-1}{2}-n} \mathlarger{\mathlarger{\mathlarger{\sl}}} \leq \\
\mathlarger{\mathlarger{\mathlarger{\sl}}}
\mathlarger{\sum}^{\frac{m(k-1)-1}{2}}_{n=2}\frac{(2\pi)^n}{n!}\(\prod^{r-1}_{j=1}\frac{\Gamma((m-j)(k-1)-n)}{\Gamma((m-j)(k-1))}\)N^{-\frac{n}{2}-1}z^{\frac{m(k-1)-1}{2}-n} \mathlarger{\mathlarger{\mathlarger{\sl}}} < \\
e^{2 \pi} \cdot \prod^{r-1}_{j=1}\frac{1}{((m-j)(k-1)-1)((m-j)(k-1)-2)}\cdot N^{-5/2}
\end{multline}
on the unit circle. In particular, we see that 
\begin{equation}\label{useRouche}
|H_{m,f}(z)-M_{m,f}(z)| < e^{2\pi}\cdot N^{-\frac{5}{2}}
\end{equation}
We also have 
\[
|M_{m,f}(z)| \geq N^{-1}-2 \pi\(\prod^{r-1}_{j=1}\frac{1}{((m-j)(k-1)-1)!}\)N^{-\frac{3}{2}} \geq N^{-1}-2 \pi N^{-\frac{3}{2}}.
\]
on the unit circle.
Thus, by comparing the highest degree terms with respect to $N$ of $M_{m,f}(z)$ and $H_{m,f}(z)-M_{m,f}(z)$ with respect to $N$ using \eqref{truncate} and \eqref{useRouche},
we see that
when $k+N$ is large enough we have
\begin{equation}\label{useRouche2}
|H_{m,f}(z)-M_{m,f}(z)|<|M_{m,f}(z)|
\end{equation}
on the unit circle.
 Thus, we may apply 
Rouch\'e's theorem to see that $H_{m,f}(z)$ and $M_{m,f}(z)$ have the same number of zeros in
 $\{z \in \C \ : \ |z|<1\}$. Since $M_{m,f}(z)$ has a zero of multiplicity $\frac{m(k-1)-1}{2}-1$ at $z=0$ and a simple zero in $\{z \in \C \ : \ |z|<1\}$ when $k$ is large enough, all of the zeros of $H_{m,f}(z)$ satisfy $|z|<1$ in this case.
\end{proof}

We now prove the analogue of Lemma~\ref{preplemma} for $Q_{m,f}(z)$, which will use Rouch\'e's theorem in a similar fashion.
\begin{lemma}\label{mainlemma}
Suppose that $m>1$ is an odd integer with $m=2r-1$, that $k$ and $N$ are positive integers with $N$ squarefree and $k \geq 6$, and that $f \in S_k(N)$ is a newform. 
If $k$ is large enough, then $Q_{m,f}(z)$ has all of its zeros in $\{z \in \C \ : \ |z| < 1\}$. If $k>2(\operatorname{log}_2(13e^{2\pi}/9)+m)+1$, then the result holds if $N$ is large enough (depending on $m$ and $k$). 
\end{lemma}
\begin{proof}
If $|z|=1$, by the definitions of \eqref{Qfm} and \eqref{eq:Hmf} we have
\begin{multline}
|H_{m,f}(z)-Q_{m,f}(z)| \leq \\
 \(\mathlarger{\sum}^{\mathsmaller{\frac{m(k-1)-1}{2}-1}}_{n=0}N^{-\frac{n}{2}-1}|L_{m,f}(m(k-1)-n)-1|\frac{(2 \pi)^n}{n!}\prod^{r-1}_{j=1}\frac{\Gamma((m-j)(k-1)-n)}{\Gamma((m-j)(k-1))}\)+\\
 \frac{1}{2}N^{-\frac{m^2(k-1)}{2}-1}\Big{|}L_{m,f}\(\frac{m(k-1)}{2}+1\)-1\Big{|}\frac{(2 \pi)^{\frac{m(k-1)-1}{2}}}{(\frac{m(k-1)-1}{2})!}\mathlarger{\prod}^{r-1}_{j=1}\frac{\Gamma\(\frac{m(k-1)-1}{2}-j(k-1)\)}{\Gamma((m-j)(k-1))}.
\end{multline}
By Lemma~\ref{lastprep}, we have
\[
|L_{m,f}(m(k-1)-n)-1| < \frac{13}{9} \times 2^{m-\frac{k-1}{2}}
\]
for $n \in \{0,\dots,\frac{k(m-1)-2m}{2}\}$, and 
by Lemma~\ref{convex}, there exists a constant $D_m$ such that
\[
 |L_{m,f}(m(k-1)-n)-1| \leq D_m \(\operatorname{log}^{m+1}\(k+\Big|\(1+\frac{m(k-1)}{2}\)\Big|+2\)\)
 \]
 for $n \in \{\frac{k(m-1)-2m}{2}+1,\dots,\frac{(m(k-1)-1)}{2}-1\}$. Thus, we have
 \begin{multline}\label{multline}
|H_{m,f}(z)-Q_{m,f}(z)| \leq \\
\frac{13}{9} \times 2^{m-\frac{k-1}{2}}e^{2\pi}N^{-1}+
\(D_m\operatorname{log}^{m+1}\(k+\Big|\(1+\frac{m(k-1)}{2}\)\Big|+2\)+1\)e^{2 \pi}N^{-\frac{k(m-1)+2m}{2}-1}\times \\
\mathlarger{\prod}^{r-1}_{j=1}\frac{\Gamma\(\frac{(m+1)(k-1)}{2}-j(k-1)\)}{\Gamma((m-j)(k-1))}.
\end{multline}
The proof of Lemma~\ref{preplemma} then shows that
\begin{equation}\label{useRouche3}
|H_{m,f}(z)-Q_{m,f}(z)|<|H_{m,f}(z)|
\end{equation}
on the unit circle when $k$ is large enough.
The result follows for large enough $k$ from applying Rouch\'e's theorem and Lemma~\ref{preplemma}.

We now consider what happens when $N$ is large. Since $k>2(\operatorname{log}_2(13e^{2\pi}/9)+m)+1$,  by applying \eqref{multline} and comparing the highest order terms of $H_{m,f}(z)$ and $H_{m,f}(z)-Q_{m,f}(z)$ with respect to $N$ (in fact $N^{-1}$), we see that \eqref{useRouche3} holds for large enough $N$ on the unit circle. Thus, our result follows in this case by arguing as above using Lemma~\ref{preplemma} and Rouch\'e's theorem.
\end{proof}
We now prove Theorem~\ref{main theorem} when $m$ is odd.
\begin{proof}[Proof of Theorem~\ref{main theorem} for odd symmetric powers]
Suppose that $k$ and $N$ are positive integers with $N$ squarefree, that $m > 1$ is an odd integer, and that $f \in S_{k}(N)$ is a newform. Let $R_{m,f}(z)$, $P_{m,f}(z)$ and $Q_{m,f}(z)$ be defined as in \eqref{symmperpoly}, \eqref{realcoeffpoly} and \eqref{Qfm}, respectively. It suffices to show that the zeros of $P_{m,f}(z)$ all lie on the unit circle when either $k$ is large enough, or $N$ is large enough assuming that $k>2(\operatorname{log}_2(13e^{2\pi}/9)+m)+1$.
In either case, we see by Lemma~\ref{mainlemma} that the zeros of $Q_{m,f}(z)$ are in $\{z \in \C \ : \ |z|<1\}$. Theorem~\ref{main theorem} then follows from Theorem~\ref{ElG,R} and \eqref{useELG,R} using $\lambda=\epsilon_{m,f}$.
\end{proof}
\section{Proof of Theorem~\ref{main theorem} for even symmetric powers}
In this section, we sketch the proof of Theorem~\ref{main theorem} when $m$ is even. 
Suppose that $m=2r$ and recall the definition of $C_{m,f}$ from Section $3$ (after \eqref{perpolyfunc}). 
In this section, we replace $C_{m,f}$ with
\[
c_{m,f}:=\pi^{\frac{m(k-1)}{2}}
\Gamma\(\frac{m(k-1)}{2}-\floor{\frac{r(k-1)}{2}}\)^{-1}C_{m,f}
\]
We then define $p_{m,f}(z):=c_{m,f}R_{m,f}(z/i\sqrt{N})$. As in \eqref{realcoefffunc}, we have
\[
p_{m,f}(1/z)=\epsilon_{m,f}z^{1-m(k-1)}p_{m,f}(z),
\]
and we define 
\[
q_{m,f}(z):=
c_{m,f} \mathlarger{\sum}^{\mathsmaller{\frac{m(k-1)-2}{2}}}_{n=0} \binom{m(k-1)-1}{n} L_{m,f}^{*}\(m(k-1)-n\)z^{\frac{m(k-1)-2}{2}-n}.
\]
Similarly as in \eqref{useELG,R}, we have
\[
p_{m,f}(z)=z^{\frac{m(k-1)-1}{2}}q_{m,f}(z)+z^{\frac{m(k-1)-2}{2}}q_{m,f}(1/z).
\]
By defining
\begin{multline}
h_{m,f}(z):=\\
\sum^{\frac{m(k-1)-2}{2}}_{n=0} \pi^{-\frac{m(k-1)-n}{2}} 
\Gamma\(\mfrac{m(k-1)-n}{2}-\floor{\mfrac{r(k-1)}{2}}\)
\mfrac{(2\pi)^n}{n!}\prod^{r-1}_{j=1}\mfrac{\Gamma((m-j)(k-1)-n)}{\Gamma((m-j)(k-1))}z^{\frac{m(k-1)-1}{2}-n}.
\end{multline}
and arguing as in the proof of Lemma~\ref{preplemma}, we obtain the following result.

\begin{lemma}
If $k+N$ is sufficiently large (independent of $m$) then all of the zeros of $h_{m,f}(z)$ satisfy $|z|<1$.
\end{lemma}

We then argue as in the proof of Lemma~\ref{mainlemma} to obtain the following lemma.
\begin{lemma} 
If $k$ is large enough, then $q_{m,f}(z)$ has all of its zeros in $\{z \in \C \ : \ |z| < 1\}$. If $k>2(\operatorname{log}_2(13e^{2\pi}/9)+m)+1$, then the result holds if $N$ is large enough (depending on $m$ and $k$).
\end{lemma}
The proof of Theorem~\ref{main theorem} for even symmetric powers then proceeds as in Section $3$ with $H_{m,f}(z)$ and $Q_{m,f}(z)$ replaced by $h_{m,f}(z)$ and $q_{m,f}(z)$, respectively.

\section{Acknowledgements}
The authors thank Jesse Thorner for helpful conversations.
 \bibliographystyle{amsalpha}
 \bibliography{main} 
 \end{document}